\documentclass{aims}

\usepackage{amsmath}
\usepackage{paralist}
\usepackage{mathrsfs}

\def\publname{~}
\def\volinfo{~}
\def\pageinfo{~}
\usepackage[colorlinks=true,unicode]{hyperref}
\hypersetup{urlcolor=blue, citecolor=red}

\def\currentvolume{X}
 \def\currentissue{X}
  \def\currentyear{200X}
   \def\currentmonth{XX}
    \def\ppages{X--XX}

\newtheorem{theorem}{Theorem}[section]

\newtheorem{lemma}[theorem]{Lemma}

\theoremstyle{definition}

\newtheorem{remark}{Remark}

\newtheorem*{example}{Example}

\title[Polynomial criteria for \textit{v-}sufficiency of map-germs]%
{Polynomial reformulation of the Kuo criteria for \textit{v-}sufficiency
of map-germs}

\author[Victor Kozyakin]{}


 \email{kozyakin@iitp.ru}

\thanks{The author is supported by the Russian Foundation for Basic
Research, project no. 10-01-93112.}

\newcommand{\grad}{\mathop{\mathrm{grad}}\nolimits}

\begin{document}
\maketitle

\centerline{\scshape Victor Kozyakin}
\medskip
{\footnotesize
 \centerline{Institute for Information Transmission Problems}
 \centerline{Russian Academy of Sciences}
  \centerline{Bolshoj Karetny lane 19, Moscow 127994 GSP-4, Russia}
} 

\bigskip

\begin{center}
Dedicated to Peter Kloeden on his 60th birthday\\
\bigskip

\end{center}

\begin{abstract}
In the paper a set of necessary and sufficient conditions for
\textit{v-}sufficiency (equiv. \textit{sv-}sufficiency) of jets
of map-germs $f:(\mathbb{R}^{n},0)\to (\mathbb{R}^{m},0)$ is
proved which generalize both the Kuiper-Kuo and the Thom
conditions in the function case ($m=1$) so as the Kuo
conditions in the general map case ($m>1$). Contrary to the Kuo
conditions the conditions proved in the paper do not require to
verify any inequalities in a so-called horn-neighborhood of the
(a'priori unknown) set $f^{-1}(0)$. Instead, the proposed
conditions reduce the problem on \textit{v-}sufficiency of jets
to evaluating the local {\L}ojasiewicz exponents for some
constructively built polynomial functions.
\end{abstract}

\section{Introduction}

In theory of dynamical systems and nonlinear analysis quite a
number of problems depending on parameters require analyzing
the structure of the set of solutions of nonlinear equations,
the number of variables in which exceeds the number of
equations. As a rule, arising equations are rather complicated
for investigation and need to be simplified in one or another
way. Clearly, such ``simplification'' may lead as to correct
conclusions about the structure of solutions as to wrong ones.
Often, small solutions of equations are of interest. In this
case one of the most popular methods of simplification of
equations is their truncation, when one casts out high order
terms in power-series expansions of the corresponding
equations. In the paper polynomial necessary and sufficient
conditions are proved allowing to judge in which cases, in
analysis of systems of real nonlinear equations of finite
smoothness, truncation is permissible.

Given a map $f:\mathbb{R}^{n}\to \mathbb{R}^{m}$ with $f(0)=0$,
let us consider the set of solutions of the equation
\begin{equation}\label{Eq11}
f(x)=0.
\end{equation}
Even locally, this set is very complicated in general. As
usual, the map $f$ is called $C^{k}$-smooth if all its
components have continuous partial derivatives up to the order
$k$ inclusive. If $f\in C^{k}$ in a neighborhood of the origin
then for each $r\leq k$ it is defined the $r$-th Taylor
polynomial $f^{(r)}(x)$ of $f(x)$ about the point $x=0$ which
will be called the \emph{$r$-truncation} of $f(x)$. Transition
from equation \eqref{Eq11} to the truncated equation
\begin{equation}\label{Eq13}
f^{(r)}(x)=0
\end{equation}
is similar to the first-approximation method in the theory of
stability and to the method of studying bifurcations by the
passage to linearized equations in nonlinear analysis, and so
on. As demonstrates the next example the sets of solutions of
equations \eqref{Eq11} and \eqref{Eq13} may be topologically
different.

\begin{example}\label{ExTrunk}
Let us discard in the next equations
$$
x_{1}^{2}-2x_{1}x_{2}^{2}+x_{1}^{4}+x_{2}^{4}+x_{2}^{8}=0,\qquad
x_{1}^{2}-2x_{1}x_{2}^{2}+x_{1}^{4}+x_{2}^{4}-x_{2}^{8}=0
$$
the terms of order higher than $4$, that is perform the
$4$-truncation of the left-hand parts. Then the truncated
equation
$$
x_{1}^{2}-2x_{1}x_{2}^{2}+x_{1}^{4}+x_{2}^{4}=
(x_{1}-x_{2}^{2})^{2}+x_{1}^{4}=0
$$
has a single solution, $x_{1}=x_{2}=0$. The first of the full
equations also has the same single solution, $x_{1}=x_{2}=0$,
while the second of the full equations has a continuum of
solutions, $x_{1}=x_{2}^{2}$. Thus, truncation of equations is
not always permissible.
\end{example}

Therefore it is natural to ask when the structure of the
zero-set of the truncated map $f^{(r)}$ is similar to that of
the full map $f$. This problem concerns the property of
sufficiency of jets. Roughly speaking, sufficiency of jets is
the property that all maps with the same truncation have the
same structure.

Following to \cite{BekkaKoike:97} we recall briefly some
definitions and results on sufficiency of jets. Let
$\mathscr{E}_{[k]}(n,m)$ denote the set of $C^{k}$ map-germs
$f:(\mathbb{R}^{n}, 0)\to (\mathbb{R}^{m}, 0)$. Given $r\leq
k$, let $j^{r}f(0)$ denote the $r$-jet of $f\in
\mathscr{E}_{[k]}(n,m)$ at $0\in \mathbb{R}^{n}$ which can be
identified with the polynomial $f^{(r)}$, and let $J^{r}(n,m)$
denote the set of $r$-jets in $\mathscr{E}_{[k]}(n,m)$. We say
$f,g\in \mathscr{E}_{[k]}(n,m)$ are $C^{0}$-equivalent, if
there is a local homeomorphism $h: (\mathbb{R}^{n}, 0)\to
(\mathbb{R}^{n}, 0)$ such that $f=g\circ h$. We further say
$f,g\in \mathscr{E}_{[k]}(n,m)$ are \textit{v-}equivalent
(resp. \textit{sv-}equivalent), if $f^{-1}(0)$ is homeomorphic
to $g^{-1}(0)$ as germs at $0\in \mathbb{R}^{n}$ (resp. there
is a local homeomorphism $h : (\mathbb{R}^{n},0)\to
(\mathbb{R}^{n},0)$ such that $h(f^{-1}(0))=g^{-1}(0))$. Given
$r\leq k$, we call an $r$-jet $w\in J^{r}(n,m)$
$C^{0}$-sufficient (resp. \textit{v-}sufficient,
\textit{sv-}sufficient) in $\mathscr{E}_{[k]}(n,m)$, if any two
maps $f,g\in \mathscr{E}_{[k]}(n,m)$ with
$j^{r}f(0)=j^{r}g(0)=w$ are $C^{0}$-equivalent (resp.
\textit{v-}equivalent, \textit{sv-}equivalent).

Clearly, $C^{0}$-sufficiency of jets implies
\textit{sv-}sufficiency, while the latter implies
\textit{v-}sufficiency. In fact, according to D.J.A. Trotman
and L.C. Wilson \cite{TrotWil:99}, \textit{v-}sufficiency is
equivalent to \textit{sv-}sufficiency.

Concerning $C^{0}$-sufficiency of jets in the function case
(i.e. $m=1$), we have

\begin{theorem}[N. Kuiper~\cite{Kuiper72}, T.-C. Kuo~\cite{Kuo69},
J. Bochnak \& S.
{\L}ojasiewicz~\cite{BochLoj71}]\label{Th-KKBL} For $f\in
\mathscr{E}_{[r]}(n,1)$, the jet $j^{r}f(0)$ is
$C^{0}$-sufficient in $\mathscr{E}_{[r]}(n,1)$ if and only if
there are positive numbers $C, \varepsilon$ such that
\begin{equation}\label{Eq-KuiKup}
|\grad f(x)|\geq C|x|^{r-1}\quad\textrm{for}\quad |x|<\varepsilon.
\end{equation}

For $f\in \mathscr{E}_{[r+1]}(n,1)$, the jet $j^{r}f(0)$ is
$C^{0}$-sufficient in $\mathscr{E}_{[r+1]}(n,1)$ if and only if
there are numbers $C,\delta, \varepsilon>0$ such that
\begin{equation}\label{Eq-KuiKupD}
|\grad f(x)|\geq C|x|^{r-\delta}\quad\textrm{for}\quad |x|<\varepsilon.
\end{equation}
\end{theorem}

K. Bekka and S. Koike \cite{BekkaKoike:97} proved that the
Kuiper-Kuo condition \eqref{Eq-KuiKup} is equivalent to the
following Thom condition: there are numbers $K,\varepsilon>0$
such that
\begin{equation}\label{Eq-Thom}
\sum_{i<j}\left| x_{i}\frac{\partial
f}{\partial x_{j}}-x_{j}\frac{\partial f}{\partial
x_{i}}\right|^{2}+|f(x)|^{2}\geq K|x|^{2r}\quad\textrm{for}\quad |x|<\varepsilon.
\end{equation}

Verification of the Kuiper-Kuo conditions \eqref{Eq-KuiKup} and
\eqref{Eq-KuiKupD}, so as of the Thom condition
\eqref{Eq-Thom}, may be reduced to the problem on evaluation of
the rate of growth of a polynomial about one of its roots,
which is equivalent to calculation of the so-called local
{\L}ojasiewicz exponents of a polynomial. Recall, that
according to the {\L}ojasiewicz theorem
\cite{Lojas:59,Lojas:91,Malg:e} for any polynomial
$p:\mathbb{R}^{n}\to\mathbb{R}$ with $p(0)=0$ there are
constants $C,\kappa>0$ such that
$$
|p(x)|\ge C|x|^{\kappa}
$$
in a neighborhood of the zero root. The least $\kappa$ for
which the above inequality holds is called the \emph{local
{\L}ojasiewicz exponent} for $p$ and is denoted by
$\mathscr{L}_{0}(p)$. If the zero root of $p$ is isolated then
such a least value of $\kappa$ exists and is rational
\cite{Gorin:61:e,Lojas:59,Lojas:91,Malg:e}. Moreover, in this
case $\mathscr{L}_{0}(p)\le (d-1)^{n}+1$ \cite{Gwozd:CMH99}
where $d$ is the degree of $p$. There is quite a number of
publications devoted to evaluation of the {\L}ojasiewicz
exponent, see, e.g.,
\cite{Abder:05,ChadKra:APM2:97,AcuKur:APM:05,BarKPl:05,BarPl:03,Gwozd:CMH99,Kollar:99,Lenar:APM:99}
and the bibliography therein.

Concerning \textit{v-}sufficiency (equiv.
\textit{sv-}sufficiency) of jets in the general map case (i.e.
$n\ge m$ but otherwise arbitrary), we have

\begin{theorem}[T.-C. Kuo~\cite{Kuo72}]\label{Th-Kuo}
For $f=\left(f_{1},f_{2},\dots,f_{m}\right)\in
\mathscr{E}_{[r]}(n,m)$ with $n\geq m$, the jet $j^{r}f(0)$ is
\textit{v-}sufficient (equiv. \textit{sv-}sufficient) in
$\mathscr{E}_{[r]}(n,m)$ if and only if there are numbers
$C,\varepsilon,\sigma>0$ such that
\begin{equation}\label{Eq-Kuo}
\mathscr{D}(\grad f^{{(r)}}_{1}(x),\grad f^{{(r)}}_{2}(x),\dots,\grad
f^{{(r)}}_{m}(x))\geq C|x|^{r-1}
\end{equation}
in $\mathscr{H}_{r}(f^{{(r)}};\sigma)\cap\{|x|<\varepsilon\}$.

For $f=\left(f_{1},f_{2},\dots,f_{m}\right)\in
\mathscr{E}_{[r+1]}(n,m)$ with $n\geq m$, the jet $j^{r}f(0)$
is \textit{v-}sufficient (equiv. \textit{sv-}sufficient) in
$\mathscr{E}_{[r+1]}(n,m)$ if and only if for any polynomial
map $g=\left(g_{1},g_{2},\dots,g_{m}\right)$ of degree $r+1$
satisfying $j^{r}g(0)=j^{r}f(0)$ there are numbers
$C,\delta,\varepsilon,\sigma>0$, all depending on $g$, such
that
\begin{equation}\label{Eq-KuoD}
\mathscr{D}(\grad f^{{(r)}}_{1}(x),\grad f^{{(r)}}_{2}(x),\dots,\grad
f^{{(r)}}_{m}(x))\geq C|x|^{r-\delta}
\end{equation}
in $\mathscr{H}_{r+1}(g;\sigma)\cap\{|x|<\varepsilon\}$.
\end{theorem}

In the above theorem, $\mathscr{H}_{s}(f;\sigma)$ denotes the
horn-neighbourhood of $f^{-1}(0)$,
$$
\mathscr{H}_{s}(f;\sigma)=\left\{x\in
\mathbb{R}^{n}: |f(x)|<\sigma|x|^{s}\right\},
$$
and
\begin{equation}\label{E-defD}
\mathscr{D}(v_{1},\dots,v_{m})=
\min_{i}\left\{\textrm{distance of~}v_{i}\textrm{~to~}V_{i}\right\}
\end{equation}
where $V_{i}$ is the span of the $v_{j}$'s, $j\neq i$.

Unfortunately, verification of the Kuo conditions
\eqref{Eq-Kuo} and \eqref{Eq-KuoD} is not as ``simple'' as
verification of the Kuiper-Kuo conditions \eqref{Eq-KuiKup},
\eqref{Eq-KuiKupD} or the Thom condition \eqref{Eq-Thom}. The
first problem here, not the major one, is that the function
$\mathscr{D}(v_{1},\dots,v_{m})$ is not defined explicitely, by
a ``simple'' formula. This causes problems in practical
evaluation of $\mathscr{D}(v_{1},\dots,v_{m})$. The second
problem, which is more serious, is that one need evaluate the
values of $\mathscr{D}(\grad f^{{(r)}}_{1}(x),\grad
f^{{(r)}}_{2}(x),\dots,\grad f^{{(r)}}_{m}(x))$ not in a
neighborhood of the origin but in horn-neighbourhoods of the
sets $(f^{{(r)}})^{-1}(0)$ or $g^{-1}(0)$ which are a'priory
unknown in general. At last, in the case of
\textit{v-}sufficiency in $\mathscr{E}_{[r+1]}(n,m)$ one need
to verify condition \eqref{Eq-KuoD} not for a single
horn-neighbourhood but for a variety of horn-neighbourhoods
defined for infinite number of polynomial maps $g$ of degree
$r+1$ satisfying $j^{r}g(0)=j^{r}f(0)$.

Not knowing about the works of N. Kuiper, T.-C. Kuo, J. Bochnak
and S. {\L}ojasiewicz, the author had sketched in
\cite{Koz:AiT84:10:e}, and later proved in
\cite[Ch.~8]{MTS:86:e}, a bit different criteria (in a bit
different terms) for \textit{sv-}sufficiency of map-germs.

\begin{theorem}[V.S. Kozyakin {\cite{Koz:AiT84:10:e},
\cite[Ch.~8]{MTS:86:e}}]\label{Th-Koz} For $f\in
\mathscr{E}_{[r]}(n,m)$ with $n\geq m$, the jet $j^{r}f(0)$ is
\textit{sv-}sufficient in $\mathscr{E}_{[r]}(n,m)$, $r\ge2$, if
and only if there is a number $q>0$ such that
\begin{equation}\label{Eq-Koz}
|f^{(r)}(x)|^{2}|y|^{2}+|(df^{(r)})^{*}(x)y|^{2}|x|^{2}\ge
q|x|^{2r}|y|^{2}
\end{equation}
for small $x$ and all $y$.

For $f\in \mathscr{E}_{[r+1]}(n,m)$ with $n\geq m$, the jet
$j^{r}f(0)$ is \textit{sv-}sufficient in
$\mathscr{E}_{[r+1]}(n,m)$, $r\ge1$, if and only if
\begin{equation}\label{Eq-KozA}
\frac{|f^{(r)}(x)|^{2}|y|^{2}+|(df^{(r)})^{*}(x)y|^{2}|x|^{2}}%
{|x|^{2r+2}|y|^{2}}\to \infty
\end{equation}
as $x\to 0$, $x\neq0$, uniformly with respect to $y\neq 0$.
\end{theorem}

In the above theorem $(df)^{*}(x)$ denotes the matrix conjugate
to $df(x)$. Clearly, the matrix $(df)^{*}(x)$ consists of $m$
column vectors $\grad f_{j}(x)$, $j=1,2,\dots,m$. If the norm
$|\cdot|$ in the above theorem is Euclidean then all the
functions in \eqref{Eq-Koz}, \eqref{Eq-KozA} are polynomial.
Hence, to verify conditions \eqref{Eq-Koz}, \eqref{Eq-KozA} one
can apply the technique of estimating the {\L}ojasiewicz
exponent mentioned above. As can be proved by standard
reasoning \cite{Gorin:61:e,Lojas:91} condition \eqref{Eq-KozA}
is equivalent, in fact, to the following condition: there are
numbers $q,\delta>0$ such that
\begin{equation}\label{Eq-KozD}
|f^{(r)}(x)|^{2}|y|^{2}+|(df^{(r)})^{*}(x)y|^{2}|x|^{2}\ge
q |x|^{2r+2-2\delta}|y|^{2}
\end{equation}
for small $x$ and all $y$, which is similar to \eqref{Eq-KuoD}.

Remark that the technique used in proving Theorem \ref{Th-Koz}
is much the same that used in proving Theorem \ref{Th-Kuo}.
Moreover, since both theorems, Theorem \ref{Th-Kuo} and Theorem
\ref{Th-Koz}, provide necessary and sufficient conditions for
\textit{sv-}sufficiency of map-germs under the same assumptions
then condition \eqref{Eq-Kuo} must be equivalent to
\eqref{Eq-Koz} while condition \eqref{Eq-KuoD} must be
equivalent to \eqref{Eq-KozD}. Nevertheless, no direct proofs
of such an equivalence, to the best of the author's knowledge,
are known.

The aim of the present paper is quite modest. First, we would
like to reformulate the Kuo conditions \eqref{Eq-Kuo},
\eqref{Eq-KuoD} in such a way to avoid verification of any
inequalities in a horn-neighborhood of the a'priori unknown set
$f^{-1}(0)$. Second, we would like to replace the function
$\mathscr{D}(\cdot)$ in \eqref{Eq-Kuo}, \eqref{Eq-KuoD} by
something easier computable in applications.

To implement this program we firstly formulate in
Section~\ref{S-regtran} ``qualified'' versions for the notions
of regularity of the set of small non-zero solutions of
equation \eqref{Eq11} and transversality of this set to small
spheres. The corresponding notions will play a key role in the
further considerations. In Lemma~\ref{Lem31} we show also that
for polynomial maps these regularity and transversality
conditions are equivalent to each other. Then, in
Theorem~\ref{Th21} we formulate a set of equivalent to each
other conditions \eqref{Eq-KozF}, \eqref{Eq-Koz1F} for
\textit{v-}sufficiency (equiv. \textit{sv-}sufficiency) of
map-germs which are direct (and trivial) generalization of
conditions \eqref{Eq-Koz}, \eqref{Eq-KozD} from
Theorem~\ref{Th-Koz}. Here we demonstrate also that these
conditions may be treated as a natural generalization of both
the Kuo conditions \eqref{Eq-Kuo}, \eqref{Eq-KuoD} and the Thom
conditions \eqref{Eq-Thom}. At last, in Section~\ref{S-proof}
to prove Theorem~\ref{Th21} we establish equivalence between
conditions \eqref{Eq-KozF}, \eqref{Eq-Koz1F} and the Kuo
conditions \eqref{Eq-Kuo}, \eqref{Eq-KuoD}.

\section{Qualified regularity and
transversality}\label{S-regtran}

Before to start formulating main results of the paper, let us
introduce some notions.

From now on $\langle\cdot,\cdot\rangle$ stands for the
Euclidean scalar product in $\mathbb{R}^{n}$, and $|\cdot|$
denotes the corresponding norm. If
$f:\mathbb{R}^{n}\to\mathbb{R}^{m}$ is a smooth map then
$(df)^{*}(x)$ is the matrix conjugate to $df(x)$. Clearly, the
matrix $(df)^{*}(x)$ consists of $m$ column vectors $\grad
f_{j}(x)$, $j=1,2,\dots,m$.

Given a map-germ $f:\mathbb{R}^{n}\to \mathbb{R}^{m}$ with
$f(0)=0$ and an integer $p\ge 1$ let us consider two auxiliary
functions of variables $x\in\mathbb{R}^{n}$ and
$y\in\mathbb{R}^{m}$:
\begin{equation}\label{Eq-defR}
\mathscr{R}_{p}(f;x,y)= |f(x)|^{p}|y|^{p}+|(df)^{*}(x)y|^{p}|x|^{p}
\end{equation}
and
\begin{equation}\label{Eq-defT}
\mathscr{T}_{p}(f;x,y)= |f(x)|^{p}|y|^{p}+
 |(df)^{*}(x)y|^{p}|x|^{p}-|\langle (df)^{*}(x)y,x\rangle|^{p}.
\end{equation}

Note, that both the functions $\mathscr{R}_{p}(f;x,y)$ and
$\mathscr{T}_{p}(f;x,y)$ are homogeneous in $y$. These
functions are polynomials in $x$ and $y$ if $f$ is a polynomial
and $p$ is even.

Positivity of the function $\mathscr{R}_{p}(f;x,y)$ for
$y\neq0$ and small $x\neq0$ means that $|(df)^{*}(x)y|>0$ for
each $y\neq0$ and all small non-vanishing solutions $x$ of
equation \eqref{Eq11}, that is the derivative of the map $f(x)$
is regular on small solutions $x$ of equation \eqref{Eq11}. So,
the inequality $\mathscr{R}(f;x,y)>0$ for $x,y\neq0$ may be
treated as a condition of regularity \cite{RoFuks:e} of small
non-zero solutions of equation \eqref{Eq11}. Therefore the
relation
\begin{equation}\label{Eq-defreg}
\mathscr{R}_{p}(f;x,y)\ge C|x|^{pq}|y|^{p},
\end{equation}
valid for small $x$ and all $y$ with some $C,q>0$, can be
called the \emph{condition of qualified regularity} of small
non-zero solutions of equation \eqref{Eq11}.

Similarly, positivity of the function $\mathscr{T}_{p}(f;x,y)$
for $y\neq0$ and small $x\neq0$ means that
$|(df)^{*}(x)y|\cdot|x|> |\langle(df)^{*}(x)y,x\rangle|$ for
each $y\neq0$ and all small non-vanishing solutions $x$ of
equation \eqref{Eq11}. The latter inequality is an algebraic
representation of the fact that the set of small solutions of
equation \eqref{Eq11} is transversal to any small sphere
$|x|=\varepsilon$ \cite{RoFuks:e}. Therefore the relation
\begin{equation}\label{Eq-deftrans}
\mathscr{T}_{p}(f;x,y)\ge C|x|^{pq}|y|^{p},
\end{equation}
valid for small $x$ and all $y$ with some $C,q>0$, can be
called the \emph{condition of qualified transversality} of
small non-zero solutions of equation \eqref{Eq11} to the small
spheres $|x|=\varepsilon$.

If the map $f$ is polynomial then the functions
$\mathscr{R}_{p}(f;x,y)$ and $\mathscr{T}_{p}(f;x,y)$ are
comparable, in a natural sense, for small $x$.

\begin{lemma}\label{Lem31}
If a map $f:\mathbb{R}^{n}\to\mathbb{R}^{m}$ with $f(0)=0$ is
polynomial then for any $p\in\mathbb{N}$ (i.e. $p$ is a natural
number) there is a constant $\mu_{p}>0$ such that
\begin{align}\label{Eq-compR}
2^{1-p}\mathscr{R}_{1}^{p}(f;x,y)&\le
\mathscr{R}_{p}(f;x,y)\le
2\mathscr{R}_{1}^{p}(f;x,y),\\
\label{Eq-compRT}
\mu_{p}\mathscr{R}_{p}(f;x,y)&\le
\mathscr{T}_{p}(f;x,y)\le\mathscr{R}_{p}(f;x,y)
\end{align}
for small $x$ and all $y$.
\end{lemma}

If a map $f$ is polynomial then by Lemma~\ref{Lem31} the set of
small non-zero solutions of equation \eqref{Eq11} is regular if
and only if it is transversal to the small spheres
$|x|=\varepsilon$, which is a well known fact \cite{Milnor:e}.
In this case the set of small non-zero solutions of equation
\eqref{Eq11} is also qualifiedly regular (with some parameter
$q>0$) if and only if it is qualifiedly transversal (with the
same parameter $q$) to the small spheres $|x|=\varepsilon$.
Moreover, all conditions \eqref{Eq-defreg} and
\eqref{Eq-deftrans} with a given $q>0$ but different
$p\in\mathbb{N}$ are equivalent to each other.

\section{Main results}\label{S-main}

\begin{theorem}\label{Th21}
For $f\in \mathscr{E}_{[r]}(n,m)$ with $n\geq m$, the jet
$j^{r}f(0)$ is \textit{v-}sufficient (equiv.
\textit{sv-}sufficient) in $\mathscr{E}_{[r]}(n,m)$ if and only
if for any $p\in\mathbb{N}$ there is a number $q>0$ such that
\begin{equation}\label{Eq-KozF}
\mathscr{K}(f^{(r)};x,y)\ge q|x|^{pr}|y|^{p}
\end{equation}
for small $x$ and all $y$ where $\mathscr{K}$ is any one of the
functions $\mathscr{R}_{p}$ or $\mathscr{T}_{p}$.

For $f\in \mathscr{E}_{[r+1]}(n,m)$ with $n\geq m$, the jet
$j^{r}f(0)$ is \textit{v-}sufficient (equiv.
\textit{sv-}sufficient) in $\mathscr{E}_{[r+1]}(n,m)$ if and
only if for any $p\in\mathbb{N}$
\begin{equation}\label{Eq-Koz1F}
\frac{\mathscr{K}(f^{(r)};x,y)}{|x|^{pr+p}|y|^{p}}\to\infty
\end{equation}
as $x\to 0$, $x\neq0$, uniformly with respect to $y\neq 0$,
where $\mathscr{K}$ is any one of the functions
$\mathscr{R}_{p}$ or $\mathscr{T}_{p}$.
\end{theorem}

Clearly, each function $\mathscr{K}(f^{(r)};x,y)$ in
Theorem~\ref{Th21} is a polynomial in $x$ and $y$, homogeneous
in $y$. This allows to simplify the formulation of Theorem
\ref{Th21} in the function case ($m=1$). Set
$$
\mathscr{R}_{p}^{*}(f^{(r)};x)= (f^{(r)}(x))^{p}+|\grad f^{(r)}(x)|^{p}|x|^{p}
$$
and
$$
\mathscr{T}_{p}^{*}(f^{(r)};x)= (f^{(r)}(x))^{p}+|\grad f^{(r)}(x)|^{p}|x|^{p}-
|\langle \grad f^{(r)}(x),x\rangle|^{p}.
$$

\begin{theorem}\label{Th22}
For $f\in \mathscr{E}_{[r]}(n,1)$, the jet $j^{r}f(0)$ is
\textit{v-}sufficient (equiv. \textit{sv-}sufficient) in
$\mathscr{E}_{[r]}(n,1)$ if and only if for any
$p\in\mathbb{N}$ there is a number $q>0$ such that
\begin{equation}\label{Eq-KozF-1}
\mathscr{K}^{*}(f^{(r)};x)\ge q|x|^{pr}
\end{equation}
for small $x$ where $\mathscr{K}^{*}$ is any one of the
functions $\mathscr{R}_{p}^{*}$ or $\mathscr{T}_{p}^{*}$.

For $f\in \mathscr{E}_{[r+1]}(n,1)$, the jet $j^{r}f(0)$ is
\textit{v-}sufficient (equiv. \textit{sv-}sufficient) in
$\mathscr{E}_{[r+1]}(n,1)$ if and only if for any
$p\in\mathbb{N}$
\begin{equation}\label{Eq-Koz1F-1}
\frac{\mathscr{K}^{*}(f^{(r)};x)}{|x|^{pr+p}}\to\infty
\end{equation}
as $x\to 0$, $x\neq0$, where $\mathscr{K}^{*}$ is any one of
the functions $\mathscr{R}_{p}^{*}$ or $\mathscr{T}_{p}^{*}$.
\end{theorem}

\begin{remark}\label{Rem-KuiKuo}
Given an analytic map-germ $h:\mathbb{R}^{n}\to\mathbb{R}^{1}$
with $h(0)=0$ and $0<\theta<1$, then the following
\emph{Bochnak-{\L}ojasiewicz inequality}
$$
|\grad h(x)|\cdot|x|\ge \theta |h(x)|
$$
holds for small $x$ \cite[Lem.~2]{BochLoj71}. Hence  for the
polynomial $f^{(r)}(x)$ in Theorem \ref{Th22} there is a number
$\gamma>0$ such that
$$
\gamma\mathscr{R}_{1}^{*}(f^{(r)};x) \le |\grad f^{(r)}(x)|\cdot|x|\le
\mathscr{R}_{1}^{*}(f^{(r)};x)
$$
for small $x$. The latter inequalities mean that conditions
\eqref{Eq-KozF-1} and \eqref{Eq-Koz1F-1} with
$\mathscr{K}^{*}=\mathscr{R}_{1}^{*}$ are equivalent to the
Kuiper-Kuo conditions \eqref{Eq-KuiKup} and \eqref{Eq-KuiKupD},
respectively.

So, conditions \eqref{Eq-KozF} and \eqref{Eq-Koz1F} in Theorem
\ref{Th21} may be treated as a natural generalization of the
Kuiper-Kuo conditions \eqref{Eq-KuiKup} and \eqref{Eq-KuiKupD},
respectively.
\end{remark}

\begin{remark}\label{Rem-Thom}
Direct verification shows that
$$
\mathscr{T}_{2}^{*}(f^{(r)};x)=
\sum_{i<j}\left| x_{i}\frac{\partial
f^{(r)}}{\partial x_{j}}-x_{j}\frac{\partial f^{(r)}}{\partial
x_{i}}\right|^{2}+|f^{(r)}(x)|^{2},
$$
and condition \eqref{Eq-KozF-1} with
$\mathscr{K}^{*}=\mathscr{T}_{2}^{*}$ is nothing else than the
Thom condition \eqref{Eq-Thom} for the map $f^{(r)}$.

So, conditions \eqref{Eq-KozF} and \eqref{Eq-Koz1F} in Theorem
\ref{Th21} may be treated also as a natural generalization for
the map case ($m>1$) of the Thom condition \eqref{Eq-Thom}.
\end{remark}

As an example of application of the formulated above theorems,
consider the well known problem on bifurcation of small
auto-oscillations from an equilibrium in a system described by
a differential equation with a parameter.

\begin{example}\label{ExHopf}
Consider the differential equation
$$
u''+\varepsilon u' + \omega^{2}u + U(\varepsilon,u,u') = 0.
$$
Let $\varepsilon$ be a small real parameter, the function
$U(\varepsilon,u,v)$ be smooth and $U(\varepsilon,0,0)\equiv
U'_{u}(\varepsilon,0,0)\equiv U'_{v}(\varepsilon,0,0)\equiv0$.
By rescaling of time and $\varepsilon$ the equation under
consideration can take the following form
\begin{equation}\label{Eq27}
u''+\frac{\lambda}{\pi} u' + u +
\frac{1}{\omega^{2}}
U\left(\frac{\lambda\omega}{\pi},u,\omega u'\right) = 0.
\end{equation}
Let $u=u(t,\lambda,\xi,\eta)$ be the solution of equation
\eqref{Eq27} satisfying the initial conditions
$u(0,\lambda,\xi,\eta)=\xi$, $u'_{t}(0,\lambda,\xi,\eta)=\eta$.
Then the problem on existence of $T$-periodic solutions of
equation \eqref{Eq27} is equivalent
\cite{AndrVitt:e,Kras:OpTrans:e} to the problem on solvability
of the following underspecified system of nonlinear equations
\begin{equation}\label{Eq28}
u(T,\lambda,\xi,\eta)=\xi,\quad u'_{t}(T,\lambda,\xi,\eta)=\eta.
\end{equation}
The left-hand parts of the last equations can be easily
evaluated, see, e.g., \cite{Kras:OpTrans:e}. Up to the second
order terms in the variables $\tau=T-2\pi$ and
$\lambda,\xi,\eta$ they have the form
$$
u(T,\lambda,\xi,\eta)=\xi+\lambda\xi+\tau\eta+\dots\,,\quad
u'_{t}(T,\lambda,\xi,\eta)=\eta-\tau\xi+\lambda\eta+\dots\,.
$$
Hence, the $2$-truncation of \eqref{Eq28} is the system of
equations
$$
\lambda\xi+\tau\eta=0,\quad \tau\xi-\lambda\eta=0.
$$
The set of solutions for these equations consists of a pair of
two-dimensional planes in the space of four-tuples
$\{\tau,\lambda,\xi,\eta\}$ having the only common point, the
zero point. One of these planes is specified by the equalities
$\tau=\lambda=0$ while the other is specified by the equalities
$\xi=\eta=0$.

Now, denote the vector $\{\tau,\lambda,\xi,\eta\}$ by $x$,
introduce an auxiliary vector $y = \{y_{1}, y_{2}\}$ and set
$$
f(x):=
\left\{u(T,\lambda,\xi,\eta)-\xi,~ u'_{t}(T,\lambda,\xi,\eta)-\eta\right\}.
$$
Then $f\in \mathscr{E}_{[2]}(4,2)$ and its $2$-truncation has
the form $f^{(2)}(x)=\left\{\lambda\xi+\tau\eta,~
\tau\xi-\lambda\eta\right\}$ from which
\begin{multline*}
\mathscr{R}_{2}(f^{(2)};x,y)= \left((\lambda\xi+\tau\eta)^{2}
+(\tau\xi-\lambda\eta)^{2}\right)\left(y_{1}^{2}+y_{2}^{2}\right)+\\
\left((\lambda y_{1}+\tau y_{2})^{2}+(\tau y_{1} -\lambda
y_{2})^{2} +(\xi y_{1} - \eta y_{2})^{2}+(\eta y_{1} +\xi
y_{2})^{2}\right)
\left(\tau^{2}+\lambda^{2}+\xi^{2}+\eta^{2}\right).
\end{multline*}
After collecting terms we get
$$
\mathscr{R}_{2}(f^{(2)};x,y)=\left((\tau^{2}+\lambda^{2})(\xi^{2}+\eta^{2})
+|x|^{4}\right)|y|^{2}\ge |x|^{4}|y|^{2}.
$$
Therefore by Theorem~\ref{Th21} the jet $j^{2}f(0)$ is
\textit{sv-}sufficient in $\mathscr{E}_{[2]}(4,2)$. Then the
set of small solutions of equations \eqref{Eq28} consists of a
pair of two-dimensional planes intersecting at the point
$\tau=\lambda=\xi=\eta=0$. Existence of one of the planes of
solutions of equations \eqref{Eq28} is obvious, it is the plane
$\xi=\eta=0$ corresponding to the trivial periodic solution
$u(t)\equiv 0$ of equation \eqref{Eq27}. Existence of the
second plane of solutions of equations \eqref{Eq28}, passing
the point $\tau=\lambda=\xi=\eta=0$ but different from the
plane $\xi=\eta=0$, testifies that equations \eqref{Eq28} have
nontrivial solutions with arbitrarily small $\tau=T-2\pi$,
$\lambda$ and $\{\xi,\eta\}\neq0$. Hence, equation \eqref{Eq27}
has small nonzero periodic solutions for some arbitrarily small
values of the parameter $\lambda$, see
\cite{KozKra:DAN80:e,Kras:OpTrans:e}.
\end{example}

\section{Proofs}\label{S-proof}
Throughout this section, $O(t^{k})$ with $k\ge0$ stands for
values having an upper bound of the form $c|t|^{k}$ for small
$t$ with some $c<\infty$. Analogously, $o(t^{k})$ denotes
values of a higher order of smallness than $|t|^{k}$ for small
$t$.

Before starting to prove Theorem~\ref{Th21} let us prove first
Lemma~\ref{Lem31}.

\subsection{Proof of Lemma~\ref{Lem31}}
Inequalities \eqref{Eq-compR} are a straightforward consequence
of the following two-sided form of the power mean inequality
$$
\left(\frac{x+y}{2}\right)^{p} \le
\frac{x^{p}+y^{p}}{2} \le (x+y)^{p},\quad p\ge1,~ x,y\ge0.
$$
So, we need only to prove inequalities \eqref{Eq-compRT} for a
given $p\in\mathbb{N}$.

The right inequality in \eqref{Eq-compRT} is obvious.
Therefore, it remains only to prove the left inequality in
\eqref{Eq-compRT} which will be done by reductio ad absurdum.

If the left inequality in \eqref{Eq-compRT} is not valid then
there are $x_{i}\to0$ ($x_{i}\neq0$), $y_{i}\neq0$ and
$\eta_{i}>0$ such that
$\eta_{i}^{p}\mathscr{R}_{p}(f;x_{i},y_{i})>\mathscr{T}_{p}(f;x_{i},y_{i})$.
In particular, $\mathscr{R}_{p}(f;x_{i},y_{i})>0$. Since the
functions $\mathscr{T}_{p}(f;x,y)$ and $\mathscr{R}_{p}(f;x,y)$
are homogeneous in $y$ with the same power of homogeneity
$p\in\mathbb{N}$ then without loss of generality one may
suppose that $|y_{i}|=1$ and $y_{i}\to y_{*}$, $|y_{*}|=1$. Let
us write the following system of polynomial equalities and
inequalities\footnote{Relations \eqref{Eq31} are polynomial
since $p$ is integer and the norm $|\cdot|$ is Euclidean.}:
\begin{gather}\nonumber
\quad |f(x)|^{2p}|y|^{2p}=u^{2p}\quad
|(df)^{*}(x)y|^{2p}|x|^{2p}=v^{2p},\quad
\langle (df)^{*}(x)y,x\rangle^{2p}=w^{2p},\\
\label{Eq31} u^{p}+v^{p}=\varphi^{p},\quad
u^{p}+v^{p}-w^{p}=\psi^{p},\quad \eta^{p}\varphi^{p}>\psi^{p},\\
\nonumber |x|^{2}>0,\quad |y|^{2}>0,\quad \varphi>0,\quad
\psi\ge0,\quad \eta>0,\quad u\ge0,\quad v\ge0,\quad w\ge0.
\end{gather}
By the definition of the sequences $\{x_{i}\}, \{y_{i}\}$ and
$\{\eta_{i}\}$, the set determined by the the relations
\eqref{Eq31} is not empty, and the point
$x=\varphi=\psi=\eta=u=v=w=0$, $y=y_{*}$ belongs to its
closure. Hence, by the Curve Selection Lemma for semialgebraic
sets (see, e.g., \cite{Milnor:e}) there are a number
$\varepsilon>0$ and real analytic around the origin functions
$x(t),y(t),\varphi(t),\psi(t)$ and $\eta(t)$ satisfying the
conditions
$$
x(0)=\varphi(0)=\psi(0)=\eta(0)=0,\quad y(0)=y_{*}
$$
and
$$
x(t)\neq0,\quad  \varphi(t)>0,\quad \psi(t)\ge0,\quad
\eta(t)>0\quad\textrm{for}\quad 0<t<\varepsilon
$$
such that
$$
\psi^{p}(t)<\eta^{p}(t)\varphi^{p}(t)\quad\textrm{for}\quad 0<t<\varepsilon
$$
or, what is the same by
\eqref{Eq-defR} and \eqref{Eq-defT},
\begin{multline}\label{Eq32}
\varphi^{p}(t)=u^{p}(t)+v^{p}(t)=\\
|f(x(t))|^{p}|y(t)|^{p}+|(df)^{*}(x(t))y(t)|^{p}|x(t)|^{p}=
\mathscr{R}_{p}(f;x(t),y(t))
\end{multline}
and
\begin{multline}\label{Eq33}
\eta^{p}(t)\varphi^{p}(t)>\psi^{p}(t)=u^{p}(t)+v^{p}(t)-w^{p}(t)=\\
|f(x(t))|^{p}|y(t)|^{p}+|(df)^{*}(x(t))y(t)|^{p}|x(t)|^{p}-
|\langle (df)^{*}(x(t))y(t),x(t)\rangle|^{p}=\\
\mathscr{T}_{p}(f;x(t),y(t)).
\end{multline}
The latter relations imply
\begin{equation}\label{Eq34}
|f(x(t))|\cdot|y(t)|\le\eta(t)\varphi(t),
\end{equation}
from which by \eqref{Eq32}
\begin{equation}\label{Eq35}
|(df)^{*}(x(t))y(t)|\cdot|x(t)|\ge
\varphi(t)\left(1-\eta^{p}(t)\right)^{1/p}.
\end{equation}
Relations \eqref{Eq33} imply also
$$
|(df)^{*}(x(t))y(t)|^{p}|x(t)|^{p}
-|\langle (df)^{*}(x(t))y(t),x(t)\rangle|^{p}
\le\eta^{p}(t)\varphi^{p}(t).
$$
By dividing the both sides of the last inequality on
$|(df)^{*}(x(t))y(t)|^{p}|x(t)|^{p}$, we obtain by \eqref{Eq35}
\begin{equation}\label{Eq36}
0\le1-
\left(\frac{|\langle (df)^{*}(x(t))y(t),x(t)\rangle|}{|(df)^{*}(x(t))y(t)|\cdot|x(t)|}\right)^{p}
\le\frac{\eta^{p}(t)}{1-\eta^{p}(t)}.
\end{equation}
Because the functions $x(t), y(t), \varphi(t), \eta(t)$ are
real analytic for small $t$ then they can be represented in the
following form:
\begin{alignat}{2}\label{Eq37}
 x(t)&=x_{*}t^{q}+o(t^{q}),&\quad x_{*}&\neq0,~ q\ge1,\\
\label{Eq38}
 y(t)&=y_{*}+O(t),&\quad |y_{*}|&=1,\\
\label{Eq39}
 \varphi(t)&=\varphi_{*}t^{r}+o(t^{r}),&\quad \varphi_{*}&>0,~ r\ge1,\\
\label{Eq310}
 \eta(t)&=\eta_{*}t^{s}+o(t^{s}),&\quad \eta_{*}&>0,~ s\ge1.
\end{alignat}
Since $f(x)$ is a polynomial and the functions $x(t), y(t)$ are
analytic then the functions $(df)^{*}(x(t))y(t)$ and $f(x(t))$
are also analytic, and $f(x(t))\to0$ as $t\to0$. Therefore by
inequalities \eqref{Eq34}--\eqref{Eq36} there are integers
$k\ge1$, $l\ge0$ such that
\begin{alignat}{2}\label{Eq311}
 f(x(t))&=O(t^{k})&\quad k&\ge1,\\
\label{Eq312}
 (df)^{*}(x(t))y(t)&=h_{*}t^{l}+o(t^{l}),&\quad
h_{*}&\neq0,~ l\ge0.
\end{alignat}
Substituting now representations \eqref{Eq38}--\eqref{Eq311}
for the related functions in \eqref{Eq34} we get
$$
O(t^{k})\cdot|y_{*}+O(t)|\le
(\varphi_{*}t^{r}+o(t^{r}))(\eta_{*}t^{s}+o(t^{s})),
$$
from which (since $y_{*}\neq0$)
\begin{equation}\label{Eq313}
k\ge r+s.
\end{equation}
Similarly, substituting representations
\eqref{Eq38}--\eqref{Eq310} and \eqref{Eq312} for the related
functions in \eqref{Eq35} we get
$$
|h_{*}t^{l}+o(t^{l})|\cdot|x_{*}t^{q}+o(t^{q})|\ge
(\varphi_{*}t^{r}+o(t^{r}))\left(1-O(t)\right)^{1/p},
$$
from which (since $h_{*}, x_{*}, \varphi_{*}\neq0$)
\begin{equation}\label{Eq314}
r\ge l+q.
\end{equation}
At last, substituting representations \eqref{Eq37},
\eqref{Eq310} and \eqref{Eq312} for the related functions in
\eqref{Eq36} we get
$$
0\le 1-
\left(\frac{|\langle h_{*}t^{l}+o(t^{l}),x_{*}t^{q}+o(t^{q})\rangle|}%
{|h_{*}t^{l}+o(t^{l})|\cdot|x_{*}t^{q}+o(t^{q})|}\right)^{p}
\le c\left(\eta_{*}t^{s}+o(t^{s})\right)^{p}
$$
with some constant $c<\infty$ from which
$$
0\le 1-
\left(\frac{|\langle h_{*},x_{*}\rangle|}{|h_{*}|\cdot|x_{*}|}\right)^{p}+O(t)
\le O(t).
$$
Hence $|\langle h_{*},x_{*}\rangle|=|h_{*}|\cdot|x_{*}|$ and
therefore $h_{*}=\lambda x_{*}$ with some $\lambda\neq0$ (since
$h_{*}, x_{*}\neq0$), and by equality \eqref{Eq312}
\begin{equation}\label{Eq315}
    (df)^{*}(x(t))y(t)=\lambda x_{*}t^{l}+o(t^{l}).
\end{equation}

Let us evaluate now the function $z(t) = \langle f(x(t)),
y(t)\rangle$. Because
\begin{multline*}
z'(t)=\langle df(x(t))x'(t), y(t)\rangle + \langle
f(x(t),y'(t)\rangle=\\ =\langle x'(t),(df)^{*}(x(t))y(t)\rangle
+ \langle f(x(t),y'(t)\rangle,
\end{multline*}
then formulae \eqref{Eq37}, \eqref{Eq38}, \eqref{Eq311} and
\eqref{Eq315} imply the following equalities
\begin{multline*}
z'(t)=\langle px_{*}t^{q-1}+O(t^{q}),\lambda
x_{*}t^{l}+o(t^{l})\rangle+
\langle O(t^{k}),O(1)\rangle=\\
=\lambda p|x_{*}|^{2}t^{q+l-1}+O(t^{q+l})+O(t^{k}).
\end{multline*}
Here, by \eqref{Eq313} and \eqref{Eq314}, $k\ge q+l+s$.
Therefore $O(t^{k})=o(t^{q+l})$, and then
$$
z'(t)=\lambda p|x_{*}|^{2}t^{q+l-1}+O(t^{q+l}).
$$
By integrating the both sides of the last equality we get
\begin{equation}\label{Eq316}
     \langle f(x(t)), y(t)\rangle=z(t)=\int_{0}^{t}z'(s)\,ds=
     \lambda\frac{q}{q+l}|x_{*}|^{2}t^{q+l}+o(t^{q+l}).
\end{equation}

Now, the obvious relation $\langle f(x(t)), y(t)\rangle\le
|f(x(t))|\cdot|y(t)|$ and inequalities \eqref{Eq316},
\eqref{Eq38} and \eqref{Eq311} imply the estimate
$$
\lambda\frac{q}{q+l}|x_{*}|^{2}t^{q+l}+o(t^{q+l})
\le O(t^{k})\cdot|y_{*}+O(t)|.
$$
Since here $x_{*}, y_{*}\neq0$ then $k\le q+l$. On the other
hand, in view of \eqref{Eq313} and \eqref{Eq314} we have $k\ge
q+l+s\ge q+l+1$. A contradiction! Lemma~\ref{Lem31} is proved.

\subsection{Proof of Theorem~\ref{Th21}}

By Lemma~\ref{Lem31} the conditions \eqref{Eq-KozF} for
different $p\in\mathbb{N}$ and $\mathscr{K}=\mathscr{R}_{p}$ or
$\mathscr{K}=\mathscr{T}_{p}$ are equivalent to each other, and
the same is valid for the conditions \eqref{Eq-Koz1F}. So, to
prove Theorem~\ref{Th21} we need only to show that the Kuo
condition \eqref{Eq-Kuo} is equivalent to the condition
\eqref{Eq-KozF} with $\mathscr{K}=\mathscr{R}_{1}$:
\begin{equation}\label{Eq-KozF-R1}
|f^{(r)}(x)|\cdot|y|+|(df^{(r)})^{*}(x)y|\cdot|x|\ge q|x|^{r}|y|
\end{equation}
for small $x$ and all $y$, while the Kuo condition
\eqref{Eq-KuoD} is equivalent to the condition \eqref{Eq-Koz1F}
with $\mathscr{K}=\mathscr{R}_{1}$:
\begin{equation}\label{Eq-Koz1F-R1}
\frac{|f^{(r)}(x)|\cdot|y|+|(df^{(r)})^{*}(x)y|\cdot|x|}{|x|^{r+1}|y|}\to\infty
\end{equation}
as $x\to 0$, $x\neq0$, uniformly with respect to $y\neq 0$.

To prove equivalence between \eqref{Eq-Kuo} and
\eqref{Eq-KozF-R1} introduce first, for a given set of vectors
$v_{1},v_{2},\dots,v_{m}\in\mathbb{R}^{n}$, the quantity
$\widetilde{\mathscr{D}}(v_{1},v_{2},\dots,v_{m})$ as follows:
\begin{equation}\label{E-defDt}
\widetilde{\mathscr{D}}(v_{1},v_{2},\dots,v_{m})=
\min \left|\sum_{i=1}^{m}y_{i}v_{i}\right|,\quad v_{1},v_{2},\dots,v_{m}\in\mathbb{R}^{n},
\end{equation}
where the minimum is taken over all $m$-tuples of real numbers
$y_{1},y_{2},\dots,y_{m}$ satisfying
$\sum_{i=1}^{m}y_{i}^{2}=1$.

Represent now the vector $(df^{(r)})^{*}(x)y$ in
\eqref{Eq-KozF-R1} in the form
$$
(df^{(r)})^{*}(x)y\equiv \sum_{i=1}^{m}y_{i}\grad
f_{i}^{(r)}(x)
$$
where $y_{1},y_{2},\dots y_{m}$ are the components of the
vector $y$ and $f^{(r)}_{1},f^{(r)}_{2},\dots f^{(r)}_{m}$ are
the components of the map $f^{(r)}$. Then, taking the minimum
in the left-hand part of \eqref{Eq-KozF-R1} over all the
vectors $y$ satisfying $\sum_{i=1}^{m}y_{i}^{2}=1$, we obtain
that
\begin{equation}\label{Eq-tD}
\min_{y\neq0}\frac{\left|(df^{(r)})^{*}(x)y\right|}{|y|}=
\widetilde{\mathscr{D}}(\grad f^{(r)}_{1}(x),\grad f^{(r)}_{2}(x),\dots,
\grad f^{(r)}_{m}(x)).
\end{equation}
Therefore \eqref{Eq-KozF-R1}, for small $x$, is equivalent to
the condition:
$$
|f^{(r)}(x)|+
\widetilde{\mathscr{D}}(\grad f^{(r)}_{1}(x),\grad f^{(r)}_{2}(x),\dots,
\grad f^{(r)}_{m}(x))\cdot|x|\ge q|x|^{r}.
$$
Then, taking into account that
\begin{equation}\label{E-DDcomp}
\widetilde{\mathscr{D}}(v_{1},v_{2},\dots,v_{m})\le
\mathscr{D}(v_{1},v_{2},\dots,v_{m})\le
\sqrt{m} \widetilde{\mathscr{D}}(v_{1},v_{2},\dots,v_{m})
\end{equation}
where $\mathscr{D}$ is the function \eqref{E-defD}, see
\cite[p.~348]{TrotWil:99}, we may state that
\eqref{Eq-KozF-R1}, for small $x$,  is equivalent also to the
condition:
\begin{equation}\label{Eq-KozF-R1D}
|f^{(r)}(x)|+
\mathscr{D}(\grad f^{(r)}_{1}(x),\grad f^{(r)}_{2}(x),\dots,
\grad f^{(r)}_{m}(x))\cdot|x|\ge \tilde{q}|x|^{r}
\end{equation}
with an appropriate $\tilde{q}>0$.

Now, let \eqref{Eq-Kuo} be valid. Then for
$x\in\mathscr{H}_{r}(f^{{(r)}};\sigma)$, $|x|<\varepsilon$, the
first summand in the left-hand side of \eqref{Eq-KozF-R1D} is
greater than $\sigma|x|^{r}$. At the same time for
$x\not\in\mathscr{H}_{r}(f^{{(r)}};\sigma)$, $|x|<\varepsilon$,
by \eqref{Eq-Kuo} the second summand in the left-hand side of
\eqref{Eq-KozF-R1D} is greater than $C|x|^{r}$. So, for
$|x|<\varepsilon$,  \eqref{Eq-Kuo} implies \eqref{Eq-KozF-R1D}
with $\tilde{q}=\min\{\sigma, C\}$.

If \eqref{Eq-KozF-R1D} holds for $|x|<\varepsilon$ with some
$\varepsilon>0$ then clearly for
$x\in\mathscr{H}_{r}(f^{{(r)}};\frac{1}{2}\tilde{q})$
inequality \eqref{Eq-Kuo} will be valid with
$C=\frac{1}{2}\tilde{q}$. So, \eqref{Eq-KozF-R1D} implies
\eqref{Eq-Kuo} with $C=\frac{1}{2}\tilde{q}$.

Thus, conditions \eqref{Eq-Kuo} and \eqref{Eq-KozF-R1D} are
equivalent, and consequently the Kuo condition \eqref{Eq-Kuo}
is equivalent to \eqref{Eq-KozF-R1}.

The proof of equivalence between \eqref{Eq-KuoD} and
\eqref{Eq-Koz1F-R1} is a bit more complicated. First, to prove
that \eqref{Eq-KuoD} implies \eqref{Eq-Koz1F-R1} we will show
that \eqref{Eq-KuoD} is not valid provided that
\eqref{Eq-Koz1F-R1} is not valid. To do it, we will need the
following lemma the proof of which is relegated to
Section~\ref{SS-Lem42} below.

\begin{lemma}\label{Lem42}
Let the map $f^{(r)}(x)$ do not satisfy \eqref{Eq-Koz1F-R1}.
Then there are $x_{i}\to0$ ($x_{i}\neq0$), $y_{i}\to0$ and a
uniform polynomial $h:\mathbb{R}^{n}\to\mathbb{R}^{m}$ of
degree $r+1$ such that for the map $g(x)=f^{(r)}(x)+h(x)$ the
following estimates hold
\begin{equation}\label{Eq42}
    |g(x_{i})|\le c|x_{i}|^{r+1+\delta'},\quad
    |(dg)^{*}(x_{i})y_{i}|\le c|y_{i}|\cdot|x_{i}|^{r+\delta'}
\end{equation}
with some $\delta'>0$ and $c<\infty$.
\end{lemma}

Now, let $\{x_{i}\}$ be a sequence defined by
Lemma~\ref{Lem42}. Then by the first inequality \eqref{Eq42}
for any $\sigma>0$ there is an $\varepsilon>0$ such that
\begin{equation}\label{E-xihorn}
    x_{i}\in \mathscr{H}_{r+1}(g;\sigma)\cap\{|x|<\varepsilon\}
\end{equation}
for all sufficiently large indices $i$.

By Lemma~\ref{Lem42}, $f^{(r)}(x)=g(x)-h(x)$ where
$h:\mathbb{R}^{n}\to\mathbb{R}^{m}$ is a uniform polynomial  of
degree $r+1$. Then $|(dh)^{*}(x)|\le c_{1}|x|^{r}$ with some
constant $c_{1}$, and by the second inequality \eqref{Eq42}
$$
|(df^{(r)})^{*}(x_{i})y_{i}|\le c_{2}|y_{i}|\cdot|x_{i}|^{r},\quad i=1,2,\ldots~,
$$
with some constant $c_{2}$. Therefore by \eqref{E-defDt}
$$
\widetilde{\mathscr{D}}(\grad f^{{(r)}}_{1}(x_{i}),\grad f^{{(r)}}_{2}(x_{i}),\dots,\grad
f^{{(r)}}_{m}(x_{i}))\leq c_{3}|x|^{r},\quad i=1,2,\ldots~,
$$
with some constant $c_{3}$, and by \eqref{E-DDcomp}
$$
\mathscr{D}(\grad f^{{(r)}}_{1}(x_{i}),\grad f^{{(r)}}_{2}(x_{i}),\dots,\grad
f^{{(r)}}_{m}(x_{i}))\leq c_{3}|x|^{r},\quad i=1,2,\ldots~.
$$
These last inequalities imply that for any $C,\delta>0$
\begin{equation}\label{E-dfrbound}
\mathscr{D}(\grad f^{{(r)}}_{1}(x_{i}),\grad f^{{(r)}}_{2}(x_{i}),\dots,\grad
f^{{(r)}}_{m}(x_{i}))< C|x|^{r-\delta}
\end{equation}
for all sufficiently large indices $i$.

Relations \eqref{E-xihorn} and \eqref{E-dfrbound} show that for
any choice of the numbers $C,\delta,\varepsilon,\sigma>0$
condition \eqref{Eq-KuoD} is not valid for the map $f^{(r)}$ so
as for the map $g$ determined by Lemma~\ref{Lem42}.

So, we completed the proof that non-validity of
\eqref{Eq-Koz1F-R1} implies non-validity of \eqref{Eq-KuoD},
and consequently the Kuo condition \eqref{Eq-KuoD} implies
\eqref{Eq-Koz1F-R1}.

It remains only to prove that \eqref{Eq-Koz1F-R1} implies the
Kuo condition \eqref{Eq-KuoD}. To do it, we will need the
following lemma the proof of which is relegated to
Section~\ref{SS-Lem32} below.

\begin{lemma}\label{Lem32}
Let $g:\mathbb{R}^{n}\to\mathbb{R}^{m}$ be a polynomial map of
degree $r+1$ such that $j^{r}g(0)=j^{r}f(0)$ where $f^{(r)}(x)$
satisfies the condition \eqref{Eq-Koz1F-R1}. Then
\begin{equation}\label{E-gfresim}
\frac{|g(x)|\cdot|y|+|(df^{(r)})^{*}(x)y|\cdot|x|}{|x|^{r+1}|y|}\to\infty
\end{equation}
as $x\to 0$, $x\neq0$, uniformly with respect to $y\neq 0$.
\end{lemma}

Now, let condition \eqref{Eq-Koz1F-R1} be valid. Take an
arbitrary polynomial map $g:\mathbb{R}^{n}\to\mathbb{R}^{m}$ of
degree $r+1$ satisfying $j^{r}g(0)=j^{r}f(0)$. Then by
Lemma~\ref{Lem32} relation \eqref{E-gfresim} holds. In this
case, by usual argumentation (see, e.g.
\cite{Gorin:61:e,Lojas:91}) there are positive constants
$\sigma',\delta'$ and  $\varepsilon'<1$ such that
\begin{equation}\label{E-gfresim-prim}
|g(x)|\cdot|y|+|(df^{(r)})^{*}(x)y|\cdot|x|\ge \sigma'|x|^{r+1-\delta'}|y|
\end{equation}
for $x\in\mathbb{R}^{n}$, $|x|<\varepsilon'$, and all
$y\in\mathbb{R}^{m}$.

Let $x$, $|x|<\varepsilon'$, belong to the horn-neighbourhood
$\mathscr{H}_{r+1}(g;\sigma'/2)$ of $g^{-1}(0)$. Then
$$
|g(x)|<\frac{1}{2}\sigma'|x|^{r+1}\le\frac{1}{2}\sigma'|x|^{r+1-\delta'}
$$
and by \eqref{E-gfresim-prim}
$$
|(df^{(r)})^{*}(x)y|\ge
\frac{1}{2}\sigma'|x|^{r-\delta'}|y|.
$$
Hence, by
\eqref{E-defDt}, \eqref{Eq-tD},
$$
\widetilde{\mathscr{D}}(\grad f^{(r)}_{1}(x),\grad f^{(r)}_{2}(x),\dots,
\grad f^{(r)}_{m}(x))\ge \frac{1}{2}\sigma'|x|^{r-\delta'},
$$
and by \eqref{E-DDcomp},
$$
\mathscr{D}(\grad
f^{(r)}_{1}(x),\grad f^{(r)}_{2}(x),\dots, \grad
f^{(r)}_{m}(x))\ge \frac{\sigma'}{2\sqrt{m}}|x|^{r-\delta'},
$$
for
$x\in\mathscr{H}_{r+1}(g;\sigma'/2)\cap\{|x|<\varepsilon'\}$,
which is exactly the Kuo condition \eqref{Eq-KuoD}.

So, \eqref{Eq-Koz1F-R1} implies the Kuo condition
\eqref{Eq-KuoD}, and the proof of Theorem~\ref{Th21} is
completed.

\subsection{Proof of Lemma~\ref{Lem42}}\label{SS-Lem42}
Denote by $H$ the class of polynomials in $x$ of the form
$\eta(x)=\langle x,v\rangle^{p}\langle x,w\rangle^{q}u$ where
$p+q=r+1$ and $v,w\in\mathbb{R}^{n}$. Since for such
polynomials
$$
d\eta(x)z=p\langle
x,v\rangle^{p-1}\langle x,w\rangle^{q}\langle z,v\rangle u+
q\langle
x,v\rangle^{p}\langle x,w\rangle^{q-1}\langle z,w\rangle u,
$$
then, by the identity $\langle (d\eta)^{*}(x)y,z\rangle \equiv
\langle y,d\eta(x)z\rangle$, it is valid also the identity
$$
\langle (d\eta)^{*}(x)y,z\rangle\equiv p\langle
x,v\rangle^{p-1}\langle x,w\rangle^{q}\langle z,v\rangle \langle y,u\rangle+
q\langle
x,v\rangle^{p}\langle x,w\rangle^{q-1}\langle z,w\rangle \langle y,u\rangle.
$$
Therefore
$$
(d\eta)^{*}(x)y=p\langle
x,v\rangle^{p-1}\langle x,w\rangle^{q}\langle y,u\rangle v+
q\langle
x,v\rangle^{p}\langle x,w\rangle^{q-1}\langle y,u\rangle w.
$$
The last formula will be needed below in two cases:
\begin{equation}\label{Eq43}
    (d\eta)^{*}(x)y=(r+1)\langle
x,v\rangle^{r}\langle y,u\rangle v,
\end{equation}
if $\eta(x)=\langle x,v\rangle^{r+1} u$, and
\begin{equation}\label{Eq44}
    (d\eta)^{*}(x)y=r\langle
x,v\rangle^{r-1}\langle x,w\rangle\langle y,u\rangle v+
q\langle
x,v\rangle^{r}\langle x,w\rangle\langle y,u\rangle w,
\end{equation}
if $\eta(x)=\langle x,v\rangle^{r}\langle x,w\rangle u$.

First, let us construct a polynomial $\eta_{i}(x)\in H$ such
that for the map $g(x)=f^{(r)}(x)+\eta_{i}(x)$ the second
inequality \eqref{Eq42} be valid. Since $f^{(r)}$ does not
satisfy \eqref{Eq-Koz1F-R1} then by the Curve Selection Lemma
\cite{Milnor:e} there are analytic around the origin functions
\begin{alignat}{2}\label{Eq45}
 x(t)&=ut^{\alpha}+o(t^{\alpha}),&\quad
 u&\neq0,~\alpha\ge1\textrm{~is integer},\\
\label{Eq46}
 y(t)&=v+O(t),&\quad |v|&=1,
\end{alignat}
for which
\begin{equation}\label{Eq48}
|f^{(r)}(x(t))|\le c |x(t)|^{r+1},\qquad
 |(df^{(r)})^{*}(x(t))y(t)|\le c |x(t)|^{r}.
\end{equation}
Clearly, the function $(df^{(r)})^{*}(x(t))y(t)$ is also
analytic. If it is identically zero then it suffices to set
$\eta_{1}(x)\equiv0$. In the opposite case let us represent it
in the form
\begin{equation}\label{Eq49}
    (df^{(r)})^{*}(x(t))y(t)=zt^{\gamma}+o(t^{\gamma}),\quad
    z\neq0,~\gamma\ge1\textrm{~is integer}.
\end{equation}
Then relations \eqref{Eq45} and \eqref{Eq48} imply
$\gamma\ge\alpha r$. If $\gamma>\alpha r$ then the second
inequality \eqref{Eq42} holds for $\eta_{1}(x)\equiv0$ and
$\delta'=(\gamma-\alpha r)\alpha$. So, it remains only to
consider the case when
\begin{equation}\label{Eq410}
    \gamma=\alpha r
\end{equation}
Here we have two possibilities: $\langle u,z\rangle\neq0$ and
$\langle u,z\rangle=0$.

a. Let $\langle u,z\rangle\neq0$. Set
$$
\eta_{1}(x)=\rho \langle
x,z\rangle^{r+1}v
$$
where $\rho\in\mathbb{R}^{1}$. By \eqref{Eq43}
$$
(d\eta_{1})^{*}(x)y=\rho(r+1) \langle
x,z\rangle^{r}\langle y,v\rangle z,
$$
and, in view of \eqref{Eq45}, \eqref{Eq46}, \eqref{Eq49} and
\eqref{Eq410}, the following equalities hold:
\begin{multline*}
(df^{(r)})^{*}(x(t))y(t)+(d\eta_{1})^{*}(x(t))y(t)=\\
 = zt^{\alpha r}+o(t^{\alpha r})
 +\rho(r+1)\langle ut^{\alpha}+o(t^{\alpha}),z\rangle^{r}
 \langle v+O(t),v\rangle z=\\
 = zt^{\alpha r}+\rho(r+1)\langle u,z\rangle^{r}
 \langle v,v\rangle zt^{\alpha r}+o(t^{\alpha r}).
\end{multline*}
If to choose $\rho=\{(r+1)\langle u,z\rangle^{r} \langle
v,v\rangle\}^{-1}$ then for the map $f^{(r)}(x)+\eta_{1}(x)$
the second estimate \eqref{Eq42} with
$\delta'=\delta'_{1}=1/(\alpha r)$ will be valid.

b. Let $\langle u,z\rangle=0$. Set
$$
\eta_{1}(x)=\rho \langle
x,z\rangle \langle x,u\rangle^{r}v
$$
where $\rho\in\mathbb{R}^{1}$. By \eqref{Eq44}
$$
(d\eta_{1})^{*}(x)y=\rho\langle x,u\rangle^{r}\langle
y,v\rangle z + \rho r\langle x,u\rangle^{r-1}\langle
x,z\rangle \langle
y,v\rangle u,
$$
and, in view of the relations \eqref{Eq45}, \eqref{Eq46},
\eqref{Eq49} and \eqref{Eq410},
\begin{multline*}
(df^{(r)})^{*}(x(t))y(t)+(d\eta_{1})^{*}(x(t))y(t)=\\
 = zt^{\alpha r}+o(t^{\alpha r})
 +\rho\langle ut^{\alpha}+o(t^{\alpha}),u\rangle^{r}
 \langle v+O(t),v\rangle z+\\
 +\rho r \langle ut^{\alpha}+o(t^{\alpha}),u\rangle^{r-1}
 \langle ut^{\alpha}+o(t^{\alpha}),z\rangle \langle v+O(t),v\rangle u.
\end{multline*}
By supposition, the multiplier $\langle
ut^{\alpha}+o(t^{\alpha}),z\rangle$ in the last summand is of
the order $O(t^{\alpha+1})$ and therefore the whole last
summand has the order $O(t^{\alpha r+1})$. Hence
$$
(df^{(r)})^{*}(x(t))y(t)+(d\eta_{1})^{*}(x(t))y(t)=
\{1+\rho\langle u,u\rangle^{r}
 \langle v,v\rangle\}t^{\alpha r} z+O(t^{\alpha r+1}).
$$
If to choose now $\rho=-\{\langle u,u\rangle^{r} \langle
v,v\rangle\}^{-1}$ then the map $f^{(r)}(x)+\eta_{1}(x)$ and
any sequences of the elements $x_{i}=x(t_{i})$,
$y_{i}=y(t_{i})$, where $t_{i}\to0$, $t_{i}\neq0$, will satisfy
the second estimate \eqref{Eq42} with
$\delta'=\delta'_{1}=1/(\alpha r)$.

So, the map $\eta_{1}(x)$ is constructed. The map $h(x)$ will
be searched in the form $h(x)=\eta_{1}(x)+\eta_{2}(x)$, with an
$\eta_{2}(x)$ such that not to break the second inequality
\eqref{Eq42} and to satisfy simultaneously the first of these
inequalities. Denote the map $f^{(r)}(x)+\eta_{1}(x)$ by
$g_{1}(x)$. Then, by construction,
\begin{align}\label{Eq411}
  |(dg_{1})^{*}(x(t))y(t)|&\le \tilde{c} |x(t)|^{r+\delta'_{1}},\\
\label{Eq412}
 |g_{1}(x(t))|&\le \tilde{c} |x(t)|^{r+1}.
\end{align}

The function $g_{1}(x(t))$ is analytic. If it is identically
zero then it suffices to set $\eta_{2}(x)\equiv0$. In the
opposite case we let us write down the following
representations:
\begin{alignat}{2}\label{Eq413}
 g_{1}(x(t))&=wt^{\mu}+o(t^{\mu}),&\quad
 w&\neq0,~\mu\ge1\textrm{~is integer},\\
\label{Eq414}
 (dg_{1})^{*}(x(t))y(t)&=O(t^{\nu}),&\quad \nu&\ge1\textrm{~is integer}.
\end{alignat}
Relations \eqref{Eq411}, \eqref{Eq414} and \eqref{Eq45} imply
\begin{equation}\label{Eq415}
    \nu\ge \alpha r+1,
\end{equation}
while relations \eqref{Eq412}, \eqref{Eq413} and \eqref{Eq45}
imply $\mu\ge \alpha(r+1)$. If $\mu> \alpha(r+1)$ then
inequalities \eqref{Eq42} hold for $\eta_{2}(x)\equiv 0$,
$\delta'=\min\{\delta'_{1},\delta'_{2}\}$ where
$\delta'_{2}=\mu/\alpha-(r+1)$. Therefore we need only to
consider the case when
\begin{equation}\label{Eq416}
    \mu=\alpha(r+1).
\end{equation}
Let us estimate the quantity $\langle g_{1}(x(t)),y(t)\rangle$.
On the one hand, by \eqref{Eq46} and \eqref{Eq413},
\begin{equation}\label{Eq417}
    \langle g_{1}(x(t)),y(t)\rangle =
    \langle wt^{\mu}+0(t^{\mu}),v+O(t)\rangle=
    \langle w,v\rangle t^{\mu} +O(t^{\mu+1}).
\end{equation}
On the other hand,
\begin{multline*}
 \langle g_{1}(x(t)),y(t)\rangle
= \int_{0}^{t}\langle dg_{1}(x(s))x'(s),y(s)\rangle\,ds+
 \int_{0}^{t}\langle g_{1}(x(s)),y'(s)\rangle\,ds=\\
 =\int_{0}^{t}\langle x'(s),(dg_{1})^{*}(x(s))y(s)\rangle\,ds+
 \int_{0}^{t}\langle g_{1}(x(s)),y'(s)\rangle\,ds,
\end{multline*}
from which, by using power series expansions in $s$ of the
integrands and by integrating the obtained relations, we get
the equalities
\begin{multline*}
 \langle g_{1}(x(t)),y(t)\rangle =\\
= \int_{0}^{t}\langle \alpha
us^{\alpha-1}+O(s^{\alpha}),O(s^{\nu})\rangle\,ds+
 \int_{0}^{t}\langle ws^{\mu}+o(s^{\mu}),O(1)\rangle\,ds=\\
 = O(t^{\alpha+\nu})+O(t^{\alpha+\nu+1})+O(t^{\mu+1})+o(t^{\mu+1}).
\end{multline*}
By \eqref{Eq415} and \eqref{Eq416}, these last equalities imply
$\langle g_{1}(x(t)),y(t)\rangle =O(t^{\mu+1})$. Therefore, in
view of \eqref{Eq417},
\begin{equation}\label{Eq418}
    \langle w,v\rangle =0.
\end{equation}

Set now $\eta_{2}(x)=\rho\langle x,u\rangle^{r+1}w$ where
$\rho=-\langle u,u\rangle^{-(r+1)}$. Then, by \eqref{Eq45} and
\eqref{Eq413},
$$
g_{1}(x(t))+\eta_{2}(x(t))=
wt^{\alpha(r+1)}+o(t^{\alpha(r+1)})-\langle u,u\rangle^{-(r+1)}
\langle ut^{\alpha}+o(t^{\alpha}),u\rangle^{r+1}w,
$$
from which $g_{1}(x(t))+\eta_{2}(x(t))=O(t^{\alpha(r+1)+1})$.
Hence, for the map $g(x)=g_{1}(x)+\eta_{2}(x)$ and any sequence
of elements $x_{i}=x(t_{i})$ where $t_{i}\to0$, $t_{i}\neq0$,
the first estimate \eqref{Eq42} holds with
$\delta'=1/(\alpha(r+1))$.

It remains to verify validity of the second estimate
\eqref{Eq42}. By \eqref{Eq43}
$$
(d\eta_{2})^{*}(x)y=\rho(r+1)\langle
x,u\rangle^{r}\langle y,w\rangle u,
$$
and therefore (see \eqref{Eq45}, \eqref{Eq46}, \eqref{Eq414})
\begin{multline*}
 (dg)^{*}(x(t))y(t)=
 (dg_{1})^{*}(x(t))y(t)+(d\eta_{2})^{*}(x(t))y(t)=\\
 = O(t^{\nu})+\rho(r+1)\langle
ut^{\alpha}+o(t^{\alpha}),u\rangle^{r}\langle v+O(t),w\rangle
u.
\end{multline*}
Since, in view of \eqref{Eq418}, the multiplier $\langle
v+O(t),w\rangle$ in the second summand is of the order $O(t)$
then the whole second summand has the order $O(t^{\alpha
r+1})$. Then by \eqref{Eq415} $(dg)^{*}(x(t))y(t)=O(t^{\alpha
r+1})$.

So, for any sequence of pairs $\{x_{i},y_{i}\}$ where
$x_{i}=x(t_{i})$, $y_{i}=y(t_{i})$, $t_{t}\to0$, $t_{i}\neq0$,
the inequalities \eqref{Eq42} hold with
$\delta'=\min\{1/(\alpha(r+1)),1/(\alpha r)\}$. The proof of
Lemma~\ref{Lem42} is completed.

\subsection{Proof of Lemma~\ref{Lem32}}\label{SS-Lem32}
Set $\theta(x)=g(x)-f^{(r)}(x)$. Then $\theta$ is a uniform
polynomial of degree $r+1$. Therefore $|\theta(x)|\le
c|x|^{r+1}$ for sufficiently small values of $|x|$ where $c$ is
some constant. Then
$$
|g(x)|\cdot|y|+|(df^{(r)})^{*}(x)y|\cdot|x|=
|f^{(r)}(x)+\theta(x)|\cdot|y|+|(df^{(r)})^{*}(x)y|\cdot|x|,
$$
and
$$
\frac{|g(x)|\cdot|y|+|(df^{(r)})^{*}(x)y|\cdot|x|}{|x|^{r+1}|y|}\ge
\frac{|f^{(r)}(x)|\cdot|y|+|(df^{(r)})^{*}(x)y|\cdot|x|}{|x|^{r+1}|y|}-
\frac{|\theta(x)|}{|x|^{r+1}},
$$
from which
$$
\frac{|g(x)|\cdot|y|+|(df^{(r)})^{*}(x)y|\cdot|x|}{|x|^{r+1}|y|}\ge
\frac{|f^{(r)}(x)|\cdot|y|+|(df^{(r)})^{*}(x)y|\cdot|x|}{|x|^{r+1}|y|}-
c.
$$
It remains to apply formula \eqref{Eq-Koz1F-R1}.
Lemma~\ref{Lem32} is proved.


\end{document}